\documentclass[12pt, reqno]{amsart}
\usepackage{amsmath, amsthm, amscd, amsfonts, amssymb, graphicx, color}
\usepackage[bookmarksnumbered, colorlinks, plainpages]{hyperref}

\newtheorem{theorem}{Theorem}[section]
\newtheorem{lemma}[theorem]{Lemma}
\newtheorem{proposition}[theorem]{Proposition}
\newtheorem{corollary}[theorem]{Corollary}
\theoremstyle{definition}
\newtheorem{definition}[theorem]{Definition}
\newtheorem{example}[theorem]{Example}

\theoremstyle{remark}
\newtheorem{remark}[theorem]{Remark}
\numberwithin{equation}{section}

\allowdisplaybreaks
\begin{document}

\title[Characterization of duals of continuous frames in Hilbert $C^{\ast}$-modules]
{Characterization of duals of continuous frames in Hilbert $C^{\ast}$-modules}

\author[H.Ghasemi]{Hadi Ghasemi}
\address{Hadi Ghasemi \\ Department of Mathematics and Computer Sciences, Hakim Sabzevari University, Sabzevar, P.O. Box 397, IRAN}
\email{ \rm h.ghasemi@hsu.ac.ir}
\author[T.L. Shateri]{Tayebe Lal Shateri }
\address{Tayebe Lal Shateri \\ Department of Mathematics and Computer Sciences, Hakim Sabzevari University, Sabzevar, P.O. Box 397, IRAN}
\email{ \rm  t.shateri@hsu.ac.ir; shateri@ualberta.ca}
\author[A.A. Arefijamaal]{Ali Akbar Arefijamaal }
\address{Ali Akbar Arefijamaal \\ Department of Mathematics and Computer Sciences, Hakim Sabzevari University, Sabzevar, P.O. Box 397, IRAN}
\email{ \rm arefijamaal@hsu.ac.ir; arefijamaal@gmail.com}
\thanks{*The corresponding author: t.shateri@hsu.ac.ir; shateri@ualberta.ca (Tayebe Lal Shateri)}
 \maketitle

\begin{abstract}
In this paper, we investigate some characterizations of dual continuous frames and give some results about them. Also, we refer to the method of constructing a family of duals through a fixed dual and show there exists a one-to-one correspondence between duals of a continuous frame for Hilbert $C^*$-module $U$ and adjointable operators $K$ from $U$ to $L^{2}(\Omega ,\mathcal A)$. Then we check the conditions that the sum of two duals of a given continuous frame under the influence of adjointable mappings becomes a dual of it and state some results about them.
\vskip 3mm
\end{abstract}

\section{Introduction And Preliminaries}
The concept of a generalization of frames to a family indexed by some locally compact space endowed with a Radon measure was proposed by G. Kaiser \cite{KAI} and independently by Ali, Antoine and Gazeau \cite{AAG}. These frames are known as continuous frames. 

Frames for Hilbert spaces have natural generalizations in Hilbert $C^*$-modules that are generalizations of Hilbert spaces by allowing the inner product to take values in a $C^*$-algebra rather than $\mathbb C$. Note that there are many differences between Hilbert $C^*$-modules and Hilbert spaces. Unlike Hilbert space cases, not every closed submodule of a Hilbert $C^*$-module is complemented. Moreover, the well-known Riesz representation theorem for continuous functionals in
Hilbert spaces does not hold in Hilbert $C^*$-modules, which implies that not all bounded linear operators on Hilbert $C^*$-modules are adjointable.  There are many essential
differences between Hilbert space frames and Hilbert $C^*$-module frames.   In Hilbert spaces, every Riesz basis has a unique dual which is also a Riesz basis. But in Hilbert $C^*$-modules, due to
the existence of zero-divisors, not all Riesz bases have unique duals, and not every dual is a Riesz basis.\\
Frank and Larson \cite{FL} presented a general approach to the frame theory in Hilbert $C^*$-modules. We refer the readers for a discussion of frames in Hilbert $C^*$-modules to Refs. \cite{AR,JI,KAB,SH,XZ}. 

The paper is organized as follows. First, we recall the basic definition of Hilbert $C^*$-modules, and we give some properties of them which we will use in the later sections. Also, we recall the notion of continuous frames in Hilbert $C^*$-modules. In Section 2,  we describe the characterizations of the duals of a continuous frame and obtain some results about them. Finally, we give the conditions under which the sum of two duals of a continuous frame is also a dual. 

We give only a brief introduction to the theory of Hilbert $C^*$-modules to make our explanations self-contained. For comprehensive accounts, we refer to  \cite{LAN,OLS}. Throughout this paper, $\mathcal A$ shows a unital $C^*$-algebra.
\begin{definition}
A \textit{pre-Hilbert module} over unital $C^*$-algebra $\mathcal A$ is a complex vector space $U$ which is also a left $\mathcal A$-module equipped with an $\mathcal A$-valued inner product $\langle .,.\rangle :U\times U\to \mathcal A$ which is $\mathbb C$-linear and $\mathcal A$-linear in its first variable and satisfies the following conditions:\\
$(i)\; \langle f,f\rangle \geq 0$,\\
$(ii)\; \langle f,f\rangle =0$  iff $f=0$,\\
$(iii)\; \langle f,g\rangle ^*=\langle g,f\rangle ,$\\
$(iv)\; \langle af,g\rangle=a\langle f,g\rangle ,$\\
for all $f,g\in U$ and $a\in\mathcal A$.
\end{definition}
A pre-Hilbert $\mathcal A$-module $U$ is called \textit{Hilbert $\mathcal A$-module} if $U$ is complete with respect to the topology determined by the norm $\|f\|=\|\langle f,f\rangle \|^{\frac{1}{2}}$.

The \textit{Cauchy-Schwartz inequality} reconstructed in Hilbert $C^*$-modules as follow \cite{OLS}.
\begin{lemma}
(\textbf{Cauchy-Schwartz inequality}) Let $U$ be a Hilbert $C^*$-modules over a unital $C^*$-algebra $\mathcal A$. Then
\begin{equation*}
\|\langle f,g\rangle \|^{2}\leq\|\langle f,f\rangle \|\;\|\langle g,g\rangle \|
\end{equation*}
for all $f,g\in U$.
\end{lemma}
For two Hilbert $C^*$-modules $U$ and $V$ over a unital $C^*$-algebra $\mathcal A$. A map $T:U\to V$ is said to be \textit{adjointable} if there exists a map $T^{*}:V\to U$ satisfying
$$\langle Tf,g\rangle =\langle f,T^*g\rangle $$
for all $f\in U, g\in V$. Such a map $T^*$ is called the \textit{adjoint} of $T$. By $End_{\mathcal A}^*(U)$ we denote the set of all adjointable maps on $U$. It is surprising that every adjointable operator is automatically linear and bounded. We need the following result in the next sections.
\begin{theorem}\label{BL}
\cite[Theorem 2.1.4]{MT} Let $U$ and $V$ be two Hilbert $C^*$-modules over a unital $C^*$-algebra $\mathcal A$ and $T\in End_{\mathcal A}^*(U,V)$. Then The following are equivalent:\\
(i)\;$T$ is bounded and $\mathcal A$-linear.\\
(ii)\;There exists $k>0$ such that $\langle Tf,Tf\rangle\leq k\langle f,f\rangle$, for all $f\in U$.
\end{theorem}
Let $\mathcal Y$ be a Banach space, $(\mathcal X,\mu)$ a measure space, and $f:\mathcal X\to \mathcal Y$ a measurable
function. The integral of the Banach-valued function $f$ has been defined by Bochner and others. Most properties of this integral are similar to those of the integral of real-valued functions (see \cite{DAN,YOS}). Since every $C^*$-algebra and Hilbert $C^*$-module is a Banach space, hence we can use this integral in these spaces.In the sequel, we assume that $\mathcal A$ is a unital $C^*$-algebra and $U$ is a Hilbert $C^*$-module over $\mathcal A$ and $(\Omega ,\mu)$ is a measure space. Define
\begin{equation*}
L^{2}(\Omega ,A)=\lbrace\varphi :\Omega \to A\quad ;\quad \Vert\int_{\Omega}\vert(\varphi(\omega))^{*}\vert^{2} d\mu(\omega)\Vert<\infty\rbrace 
\end{equation*}
It was shown that $L^{2}(\Omega , A)$ is a Hilbert $\mathcal A$-module with the inner product
\begin{equation*}
\langle \varphi ,\psi\rangle = \int_{\Omega}\langle\varphi(\omega),\psi(\omega)\rangle d\mu(\omega),
\end{equation*}
and induced norm $\|\varphi\|=\|\langle \varphi,\varphi\rangle \|^{\frac{1}{2}}$, for any $\varphi ,\psi \in L^{2}(\Omega , A)$. \cite{LAN}

In the following, we recall the notion of continuous frames in Hilbert $C^*$-modules over a unital $C^*$-algebra $\mathcal A$, and mention some properties of these frames. For details, see \cite{GHSH1,GHSH}.
\begin{definition}
A mapping $F:\Omega \to U$ is called a continuous frame for $U$ if\\
$(i)\; F$ is weakly-measurable, i.e., the mapping $\omega\longmapsto\langle f,F(\omega)\rangle $ is measurable on $\Omega$, for any $f\in U$.\\
$(ii)$ There exist constants $A,B>0$ such that
\begin{equation}\label{eq1}
A\langle f,f\rangle \leq \int_{\Omega}\langle f,F(\omega)\rangle \langle F(\omega),f\rangle d\mu(\omega)\leq B\langle f,f\rangle  ,\quad (f\in U). 
\end{equation}
\end{definition}
The constants $A,B$ are called \textit{lower} and \textit{upper} frame bounds, respectively. The mapping $F$ is called \textit{Bessel} if the right inequality in \eqref{eq1} holds and is called \textit{tight} if  $A=B$.

For a continuous frame $F:\Omega \to U$ the following operators were defined.\\
$(i)$\;The \textit{synthesis operator} or \textit{pre-frame operator} $T_{F}:L^{2}(\Omega , A)\;\to U$ weakly defined by
\begin{equation}
\langle T_{F}\varphi ,f\rangle =\int_{\Omega}\varphi(\omega)\langle F(\omega),f\rangle d\mu(\omega),\quad (f\in U).
\end{equation}
$(ii)$\; The adjoint of $T,$ so called the \textit{analysis operator} $T^{\ast}_{F}:U\;\to L^{2}(\Omega , A)$ is given by
\begin{equation}
(T^{\ast}_{F}f)(\omega)=\langle f ,F(\omega)\rangle\quad (\omega\in \Omega).
\end{equation}
$(iii)$\; The \textit{frame operator} $S_{F}:U\;\to U$ is weakly defined by
\begin{equation}
\langle S_{F}f ,f\rangle =\int_{\Omega}\langle f,F(\omega)\rangle\langle F(\omega),f\rangle d\mu(\omega),\quad (f\in U).
\end{equation}
In \cite{GHSH} was proved the pre-frame operator $T_{F}:L^{2}(\Omega , A)\;\to U$ is well defined, surjective, adjointable $\mathcal A$-linear map and bounded with $\| T\|\leq\sqrt{B}$ . Moreover the analysis operator $T^{\ast}_{F}:U\;\to L^{2}(\Omega , A)$ is injective and has closed range. Also $S=T T^{\ast}$ is positive, adjointable, self-adjoint and invertible and $\| S\|\leq B$.

Now we recall the concept of duals of continuous frames in Hilbert $C^{\ast}$-modules \cite{GHSH3}.
\begin{definition}
Let $F:\Omega \to U$ be a continuous Bessel mapping. A continuous Bessel mapping $G:\Omega \to U$ is called a \textit{dual} for $F$ if
\begin{equation*}
f= \int_{\Omega}\langle f,G(\omega)\rangle F(\omega)d\mu(\omega),\qquad(f\in U)
\end{equation*}
or
\begin{equation}\label{eq2}
\langle f,g\rangle= \int_{\Omega}\langle f,G(\omega)\rangle\langle F(\omega),g\rangle d\mu(\omega),\qquad(f,g\in U).
\end{equation}
In this case $(F,G)$ is called a \textit{dual pair}. If $T_{F}$ and $T_{G}$ denote the synthesis operators of $F$ and $G$, respectively, then \eqref{eq2} is equivalent to $T_{F}T^{*}_{G}=I_{U}$.
\end{definition}

\begin{remark}
Let $F:\Omega \to U$ be a continuous frame for Hilbert $C^{\ast}$-module $U$. Then by reconstructin formula we have
\begin{equation*}
f=\int_{\Omega}\langle f,F(\omega)\rangle S^{-1}F(\omega) d\mu(\omega)=\int_{\Omega}\langle f,S^{-1}F(\omega)\rangle F(\omega) d\mu(\omega).
\end{equation*}
\end{remark}
Therefore, $S^{-1}F$ is a dual for $F$, which is called the \textit{canonical dual}. If a continuous frame $F$ has only one dual, it is called a \textit{Riesz-type frame}.

We need the following theorem in the next section.
\begin{theorem}
\cite[Theorem 3.4]{GHSH} Let $F:\Omega \to U$ be a continuous frame for Hilbert $C^{\ast}$-module $U$ over a unital $C^*$-algebra $\mathcal A$. Then $F$ is a Riesz-type frame if and only if the analysis operator $T^{*}_{F}:U\to L^{2}(\Omega ,A)$ is onto.
\end{theorem}      
\section{Construction of dual continuous frames}
 
In this section, we characterize the duals of a continuous frame and give some results about them. We refer to the method of constructing a family of duals through a fixed dual. We also examine the conditions so that a mapping becomes a continuous frame under the influence of a continuous frame.

In the following results, we show how we can construct a sequence of duals of a continuous frame from a given dual.
\begin{theorem}\label{f-of-Ds}
Let $F:\Omega \to U$ be a continuous frame for Hilbert $C^{\ast}$-module $U$ over a unital $C^*$-algebra $\mathcal A$ with the continuous frame operator $S_{F}$ and also assume that $G_{1}:\Omega \to U$ is a dual of $F$. Then $\lbrace G_{i}\rbrace_{i\in\mathbb{N}}$ where
\begin{equation*}
G_{i+1}(\omega)=S^{-1}_{F}F(\omega)+S_{F}G_{i}(\omega)-F(\omega)
\end{equation*}
is a sequence of duals of $F$.
\end{theorem}

\begin{proof}
Clearly $G_{i}$ is a continuous Bessel mapping for $U$, for each $i\in\mathbb{N}$. Now we use induction to continue the proof. If $i=1$, then
\begin{align*}
\int_{\Omega}\langle f,F(\omega)\rangle\langle G_{2}(\omega),g\rangle d\mu(\omega) &= \int_{\Omega}\langle f,F(\omega)\rangle \langle S^{-1}_{F}F(\omega),g\rangle d\mu(\omega)+\int_{\Omega}\langle f,F(\omega)\rangle\langle S_{F}G_{1}(\omega),g\rangle d\mu(\omega)\\
&\qquad - \int_{\Omega}\langle f,F(\omega)\rangle\langle F(\omega),g\rangle d\mu(\omega)\\
& =\langle f,g\rangle +\langle f,S_{F}g\rangle -\langle S_{F}f,g\rangle =\langle f,g\rangle,
\end{align*}
for every $f,g\in U$. Suppose that $G_{k}$ is a dual of $F$. We show that $G_{k+1}$ is also a dual of $F$.
\begin{align*}
\int_{\Omega}\langle f,F(\omega)\rangle\langle G_{k+1}(\omega),g\rangle d\mu(\omega) &= \int_{\Omega}\langle f,F(\omega)\rangle \langle S^{-1}_{F}F(\omega),g\rangle d\mu(\omega)+\int_{\Omega}\langle f,F(\omega)\rangle\langle S_{F}G_{k}(\omega),g\rangle d\mu(\omega)\\
&\qquad - \int_{\Omega}\langle f,F(\omega)\rangle\langle F(\omega),g\rangle d\mu(\omega)\\
& =\langle f,g\rangle +\langle f,S_{F}g\rangle -\langle S_{F}f,g\rangle =\langle f,g\rangle,
\end{align*}
for every $f,g\in U$.
\end{proof}
\begin{corollary}\label{seq of D}
Let $F:\Omega \to U$ be a continuous frame for Hilbert $C^{\ast}$-module $U$ over a unital $C^*$-algebra $\mathcal A$ with the continuous frame operator $S_{F}$ and let $G:\Omega \to U$ be a dual of $F$. Then $\lbrace V_{i}\rbrace_{i\in\mathbb{N}}$ is a sequence of duals of $F$, where
\begin{equation*}
V_{1}(\omega)=S^{-1}_{F}F(\omega)+S_{F}G(\omega)-F(\omega),
\end{equation*}
and
\begin{equation*}
V_{i+1}(\omega)=S^{-1}_{F}F(\omega)+S^{i+1}_{F}G(\omega)-S^{i}_{F}F(\omega),
\end{equation*}
for every $i\in\mathbb{N}$.
\end{corollary}

\begin{proof}
By Theorem \ref{f-of-Ds}, the mapping $V_{1}$ is a dual of $F$. Also
\begin{align*}
\int_{\Omega}\langle f,V_{i+1}F(\omega)\rangle F(\omega)d\mu(\omega) & =\int_{\Omega}\langle f,S^{-1}_{F}F(\omega)+S^{i+1}_{F}G(\omega)-S^{i}F(\omega)\rangle F(\omega) d\mu(\omega)\\
& = \int_{\Omega}\langle f,S^{-1}_{F}F(\omega)\rangle F(\omega)d\mu(\omega)+\int_{\Omega}\langle f,S^{i+1}_{F}G(\omega)\rangle F(\omega)d\mu(\omega)\\
&\qquad - \int_{\Omega}\langle f,S^{i}F(\omega)\rangle F(\omega)d\mu(\omega)\\
& = f+S^{i+1}_{F}f-S_{F}(S^{i}_{F}f)=f,
\end{align*}
for every $f\in U$ and $i\in\mathbb{N}$.
\end{proof}
The following theorem shows that the difference between each arbitrary dual and the canonical dual of a continuous frame can be considered as a continuous Bessel mapping.
\begin{theorem}\label{dif-duals}
Let $F:\Omega \to U$ be a continuous frame for Hilbert $C^{\ast}$-module $U$ over a unital $C^*$-algebra $\mathcal A$ with the frame operator $S_{F}$. Then a continuous Bessel mapping $G:\Omega \to U$ is a dual of $F$ if and only if there exists a continuous Bessel mapping $L:\Omega \to U$ such that
\begin{equation*}
G(\omega)=S^{-1}_{F}F(\omega)+L(\omega)\qquad(\omega\in\Omega),
\end{equation*}
where $\int_{\Omega}\langle f,F(\omega)\rangle L(\omega) d\mu(\omega) =0$ for each $f\in U$.
\end{theorem}

\begin{proof}
Let $G$ be a dual of $F$.Then $f=\int_{\Omega}\langle f,GF(\omega)\rangle F(\omega) d\mu(\omega)$ for each $f\in U$. Suppose that $L=G-S^{-1}_{F}F$. Then for each $f,g\in U$,
\begin{align*}
\int_{\Omega}\langle f,F(\omega)\rangle\langle L(\omega),g\rangle d\mu(\omega) & =\int_{\Omega}\langle f,F(\omega)\rangle\langle G-S^{-1}_{F}F(\omega),g\rangle d\mu(\omega)\\
& =\int_{\Omega}\langle f,F(\omega)\rangle\langle G(\omega),g\rangle d\mu(\omega) -\int_{\Omega}\langle f,F(\omega)\rangle\langle S^{-1}_{F}F(\omega),g\rangle d\mu(\omega)\\
& =\langle f,g\rangle -\langle f,g\rangle=0.
\end{align*}
Hence $\langle\lbrace\langle f,F(\omega)\rangle\rbrace_{\omega\in\Omega},\lbrace\langle g,L(\omega)\rangle\rbrace_{\omega\in\Omega}\rangle =0$ and $\lbrace\langle f,F(\omega)\rangle\rbrace_{\omega\in\Omega}\in R(T^{*}_{L})^{\perp} =Ker(T_{L})$, where $T_{L}$ is the pre-frame operator of $L$. This shows that $\int_{\Omega}\langle f,F(\omega)\rangle L(\omega) d\mu(\omega) =0$.

Conversly, let $G=S^{-1}_{F}F+L$ where $\int_{\Omega}\langle f,F(\omega)\rangle L(\omega) d\mu(\omega) =0$ for each $f\in U$. Then,
\begin{align*}
\int_{\Omega}\langle f,F(\omega)\rangle\langle G(\omega),g\rangle d\mu(\omega) & =\int_{\Omega}\langle f,F(\omega)\rangle\langle S^{-1}_{F}F(\omega)+L(\omega),g\rangle d\mu(\omega)\\
& =\int_{\Omega}\langle f,F(\omega)\rangle\langle S^{-1}_{F}F(\omega),g\rangle d\mu(\omega)+\int_{\Omega}\langle f,F(\omega)\rangle\langle L(\omega),g\rangle d\mu(\omega)\\
& =\langle f,g\rangle +\langle\int_{\Omega}\langle f,F(\omega)\rangle L(\omega) d\mu(\omega) ,g\rangle\\
& = \langle f,g\rangle +0 = \langle f,g\rangle,
\end{align*}
for each $f,g\in U$. Hence $G$ is a dual of $F$.
\end{proof}
\begin{example}\label{ex-seq-dual}
Assume that $\mathcal A=\Big\{\begin{pmatrix}
a&b\\
c&d
\end{pmatrix}: a,b,c,d\in \mathbb C\Big\}$ which is an unital $C^*$-algebra. We define the inner product
\begin{equation*}
\begin{array}{ll}
\langle .,.\rangle:\mathcal A\times \mathcal A\;\to \quad \mathcal A \\
\qquad\; (M,N)\longmapsto M(\overline{N})^t.
\end{array}
\end{equation*}
It is not difficult to see that $\mathcal A$ with this inner product is a Hilbert $C^*$-module over itself. Suppose that $(\Omega ,\mu)$ is a measure space where  $\Omega=[0,1]$ and $\mu$ is the Lebesgue measure. Also the mapping $F:\Omega\to \mathcal A$ is defined by 
$F(\omega)=\begin{pmatrix}
2\omega &\omega\\
\omega &3\omega
\end{pmatrix}$, for any $\omega\in \Omega$.\\
It is easy to see that $F$ is a continuous frame for $\mathcal A$ with the frame bounds $\dfrac{1}{2},\dfrac{9}{2}$ and
\begin{equation*}
S_{F}=\begin{pmatrix}
\frac{5}{3}&\frac{5}{3}\\
\frac{5}{3}&\frac{10}{3}
\end{pmatrix}, S^{-1}_{F}=\begin{pmatrix}
\frac{6}{5}&\frac{-3}{5}\\
\frac{-3}{5}&\frac{3}{5}
\end{pmatrix}, S^{-1}_{F}F(\omega)=\begin{pmatrix}
\frac{9}{5}\omega &\frac{-3}{5}\omega\\
\frac{-3}{5}\omega &\frac{6}{5}\omega
\end{pmatrix}
\end{equation*}
Also, the mapping $G:\Omega\to \mathcal A$ where
\begin{equation*}
G(\omega)=\begin{pmatrix}
\frac{33}{10}\omega -1 &\frac{9}{10}\omega -1\\
\frac{9}{10}\omega -1&\frac{27}{10}\omega -1
\end{pmatrix}
\end{equation*}
is a dual of $F$ and by Theorem \ref{seq of D}, $\lbrace V_{i}\rbrace_{i\in\mathbb{N}}$ is a sequence of duals of $F$, where
\begin{equation*}
V_{1}(\omega)=\begin{pmatrix}
\frac{102\omega -50}{15} &\frac{66\omega -50}{15}\\
\frac{69\omega -50}{10}&\frac{87\omega -50}{10}
\end{pmatrix}
\end{equation*}
and
\begin{equation*}
V_{i+1}(\omega)=\begin{pmatrix}
\frac{9}{5}\omega &\frac{-3}{5}\omega\\
\frac{-3}{5}\omega &\frac{6}{5}\omega
\end{pmatrix} +\begin{pmatrix}
\frac{5}{3}\omega &\frac{5}{3}\omega\\
\frac{5}{3}\omega &\frac{10}{3}\omega
\end{pmatrix}^{i}\begin{pmatrix}
5\omega -\frac{10}{3} &5\omega -\frac{10}{3}\\
\frac{15}{2}\omega -5 &\frac{15}{2}\omega -5
\end{pmatrix}
\end{equation*}
for all $i\in\mathbb{N}$.\\
Moreover, due to Theorem \ref{dif-duals}, each dual of $F$ is obtained as the mapping $G:\Omega\to \mathcal A$ defined by
\begin{equation*}
G_{\alpha ,\beta ,\gamma ,\delta}(\omega)=\begin{pmatrix}
(\alpha +\frac{9}{5})\omega -\frac{2}{3}\alpha &(\beta +\frac{-3}{5})\omega -\frac{2}{3}\beta\\
(\gamma +\frac{-3}{5})\omega -\frac{2}{3}\gamma &(\delta +\frac{6}{5})\omega -\frac{2}{3}\delta
\end{pmatrix}
\end{equation*}
where $\alpha ,\beta ,\gamma ,\delta\in\mathbb{R}$ are orbitrary.
\end{example}
It was shown in \cite{ZQX} that the difference between each arbitrary dual and the canonical dual of a continuous frame in Hilbert spaces can be considered as a bounded operator. In the following, we generalize it for Hilbert $C^{\ast}$-modules.
\begin{theorem}
Let $F:\Omega \to U$ be a continuous frame for Hilbert $C^{\ast}$-module $U$ over a unital $C^*$-algebra $\mathcal A$ with the frame operator $S_{F}$ and the pre-frame operator $T_{F}$ and the frame bounds $A_{F},B_{F}>0$. Then there exists a one-to-one correspondence between duals of $F$ and adjointable operators $K:U \to L^{2}(\Omega ,\mathcal A)$ such that $T_{F}K=0$.\\
In this case, $\Vert K\Vert\leq\sqrt{D}+\dfrac{1}{\sqrt{A}}$ where $D$ is the upper frame dound of $G$.
\end{theorem}

\begin{proof}
Let $G$ be a dual of $F$. For each $f\in U$ , $\omega\in\Omega$ define
\begin{align*}
K: & U \longrightarrow L^{2}(\Omega ,\mathcal A)\\
& f\longmapsto Kf
\end{align*}
where $(Kf)(\omega)=\langle f, G(\omega)-S_{F}^{-1}F(\omega)\rangle$.

It is obvious that $K$ is well-defined. Also for every $f\in U$ , $\varphi\in  L^{2}(\Omega ,\mathcal A)$
\begin{align*}
\langle Kf ,\varphi\rangle &= \int_{\Omega}\langle f,G(\omega)-S_{F}^{-1}F(\omega)\rangle(\varphi(\omega))^{*} d\mu(\omega)\\
&= \langle f,\int_{\Omega}\varphi(\omega)(G(\omega)-S_{F}^{-1}F(\omega)) d\mu(\omega)\rangle = \langle f, (T_{G}- T_{S^{-1}F})\varphi\rangle.
\end{align*}
Hence $K$ is adjointable and $K^{*}=T_{G}- T_{S^{-1}F}$. Also
\begin{align*}
(T_{F}K)(f) &= \int_{\Omega}(Kf)(\omega)F(\omega) d\mu(\omega)\\
&= \int_{\Omega}\langle f,G(\omega)-S_{F}^{-1}F(\omega)\rangle F(\omega) d\mu(\omega)\\
&= \int_{\Omega}\langle f,G(\omega)\rangle F(\omega) d\mu(\omega) -\int_{\Omega}\langle f,S_{F}^{-1}F(\omega)\rangle F(\omega) d\mu(\omega)\\
&=f-f=0.
\end{align*}
i.e., $T_{F}K=0$.

Moreover, assume that $D$ is the upper frame bound of $G$. Then for $f\in U$ and $\varphi\in  L^{2}(\Omega ,\mathcal A)$ where $\Vert\varphi\Vert\leq 1$ we have
\begin{align*}
\Vert Kf\Vert & =\sup_{\Vert\varphi\Vert\leq 1}\Vert\langle Kf,\varphi\rangle\Vert\\
&= \sup_{\Vert\varphi\Vert\leq 1}\Vert\langle f,\int_{\Omega}\varphi(\omega)(G(\omega)-S_{F}^{-1}F(\omega)) d\mu(\omega)\rangle\Vert\\
&= \sup_{\Vert\varphi\Vert\leq 1}\Vert\langle f, (T_{G}- T_{S^{-1}F})\varphi\rangle\Vert\\
& \leq\sup_{\Vert\varphi\Vert\leq 1}\Vert T_{G}- T_{S^{-1}F}\Vert\Vert\varphi\Vert\Vert f\Vert\\
& \leq(\Vert T_{G}\Vert +\Vert T_{S^{-1}F}\Vert)\Vert f\Vert =(\sqrt{D}+\dfrac{1}{\sqrt{A}})\Vert f\Vert.
\end{align*}
Conversly, let $K:U \to L^{2}(\Omega ,\mathcal A)$ be adjiontable and $T_{F}K=0$. Define
\begin{align*}
G: & \Omega \longrightarrow U\\
&\omega\longmapsto G(\omega)
\end{align*}
where $\langle f, G(\omega)\rangle=(Kf)(\omega)+\langle f,S_{F}^{-1}F(\omega)\rangle$, for each $f\in U$ , $\omega\in\Omega$

Now we show that the mapping
\begin{align*}
T_{G}: L^{2}(\Omega & ,\mathcal A) \longrightarrow U\\
&\varphi\longmapsto\int_{\Omega}\varphi(\omega)G(\omega) d\mu(\omega)
\end{align*}
is well-defined and bounded. For $\varphi\in  L^{2}(\Omega ,\mathcal A)$ and $f\in U$ which $\Vert f\Vert\leq 1$ we have
\begin{align*}
\Vert T_{G}\varphi\Vert & =\sup_{\Vert f\Vert\leq 1}\Vert\langle T_{G}\varphi ,f\rangle\Vert\\
&= \sup_{\Vert f\Vert\leq 1}\Vert\int_{\Omega}\varphi(\omega)\langle G(\omega),f\rangle d\mu(\omega)\Vert\\
&= \sup_{\Vert f\Vert\leq 1}\Vert\int_{\Omega}\varphi(\omega)(Kf)(\omega)^{*} d\mu(\omega) +\int_{\Omega}\varphi(\omega)\langle S^{-1}F(\omega),f\rangle d\mu(\omega)\Vert\\
& \leq\sup_{\Vert f\Vert\leq 1}(\Vert \langle\varphi , Kf\rangle\Vert +\Vert\langle T_{S^{-1}F}\varphi ,f\rangle\Vert)\leq(\Vert K\Vert +\Vert T_{S^{-1}F}\Vert)\Vert \varphi\Vert.
\end{align*}
This implies that $G$ is a continuous Bessel mapping. Also,
\begin{align*}
\int_{\Omega}\langle f,G(\omega)\rangle F(\omega)d\mu(\omega) &=\int_{\Omega}\langle (Kf)(\omega)\rangle F(\omega)d\mu(\omega)+\int_{\Omega}\langle S_{F}^{-1}F(\omega),f\rangle F(\omega)d\mu(\omega)\\
&=T_{F}(Kf)+f=f .
\end{align*}
Therefore $G$ is a dual of $F$.
\end{proof}
In the following, we give two interesting characterizations of the cononical dual of a continuous frame.
\begin{theorem}
Let $F:\Omega \to U$ be a continuous frame for Hilbert $C^{\ast}$-module $U$ over a unital $C^*$-algebra $\mathcal A$ with the continuous frame operator $S_{F}$ and also assume that $G:\Omega \to U$ is a dual of $F$. Then $G$ is the canonical dual of $F$ if and only if
\begin{equation*}
\langle F(\omega),G(\gamma)\rangle =\langle G(\omega),F(\gamma)\rangle,
\end{equation*}
for all $\omega ,\gamma\in\Omega$.
\end{theorem}

\begin{proof}
Let $G$ be the canonical dual of $F$. Then
\begin{align*}
\langle F(\omega),G(\gamma)\rangle  &=\langle F(\omega),S^{-1}F(\gamma)\rangle\\
&=\langle S^{-1}F(\omega),F(\gamma)\rangle =\langle G(\omega),F(\gamma)\rangle ,
\end{align*}
for all $\omega ,\gamma\in\Omega$.

Conversly, suppose that $\langle F(\omega),G(\gamma)\rangle =\langle G(\omega),F(\gamma)\rangle$, for all $\omega ,\gamma\in\Omega$. Then
\begin{align*}
\langle f,F(\gamma)\rangle  &=\langle\int_{\Omega}\langle f,F(\gamma)\rangle G(\gamma)d\mu(\omega),F(\omega)\rangle\\
&=\int_{\Omega}\langle f,F(\gamma)\rangle\langle G(\gamma),F(\omega)\rangle d\mu(\omega)\\
&=\int_{\Omega}\langle f,F(\gamma)\rangle\langle F(\gamma),G(\omega)\rangle d\mu(\omega)\\
&=\langle\int_{\Omega}\langle f,F(\gamma)\rangle F(\gamma)d\mu(\omega),G(\omega)\rangle\\
&=\langle Sf,G(\omega)\rangle =\langle f,SG(\omega)\rangle ,
\end{align*}
for every $f\in U$. This shows that $F(\omega)=SG(\omega)$. Hence $G(\omega)=S^{-1}F(\omega)$, for all $\omega\in\Omega$.

\end{proof}
\begin{theorem}
Let $F:\Omega \to U$ be a continuous frame for Hilbert $C^{\ast}$-module $U$ over a unital $C^*$-algebra $\mathcal A$ with the continuous frame operator $S_{F}$ and also assume that $G:\Omega \to U$ is a dual of $F$. Then $G$ is the canonical dual of $F$ if and only if the inequality
\begin{equation*}
\langle S_{G}f,f\rangle\leq\langle S_{D}f,f\rangle ,
\end{equation*}
holds for any dual $D$ of $F$ and every $f\in U$,  where $S_{G}$ and $S_{D}$ are the continuous frame operators of $G$ and $D$, respectively.
\end{theorem}

\begin{proof}
Assume that $G=S^{-1}_{F}F$. Since both $S^{-1}_{F}F$ and $D$ are dual frames of $F$, so for each $f\in U$
\begin{equation*}
\int_{\Omega}\langle f,D(\omega)\rangle F(\omega)d\mu(\omega)=f=\int_{\Omega}\langle f,S^{-1}_{F}F(\omega)\rangle F(\omega) d\mu(\omega),
\end{equation*}
then
\begin{equation*}
\int_{\Omega}\langle f,D(\omega)-S^{-1}_{F}F(\omega)\rangle F(\omega)d\mu(\omega)=0,
\end{equation*}
and
\begin{equation*}
\int_{\Omega}\langle f,D(\omega)-S^{-1}_{F}F(\omega)\rangle\langle S^{-1}_{F}F(\omega),f\rangle d\mu(\omega)=0.
\end{equation*}
Hence
\begin{align*}
\int_{\Omega}\langle f,D(\omega)\rangle\langle D(\omega),f\rangle d\mu(\omega) &=\int_{\Omega}\langle f,D(\omega)-S^{-1}_{F}F(\omega)+S^{-1}_{F}F(\omega)\rangle\langle D(\omega)-S^{-1}_{F}F(\omega)+S^{-1}_{F}F(\omega),f\rangle d\mu(\omega)\\
&=\int_{\Omega}\langle f,D(\omega)-S^{-1}_{F}F(\omega)\rangle\langle D(\omega)-S^{-1}_{F}F(\omega),f\rangle d\mu(\omega)\\
& \qquad +\int_{\Omega}\langle f,D(\omega)-S^{-1}_{F}F(\omega)\rangle\langle S^{-1}_{F}F(\omega),f\rangle d\mu(\omega)\\
& \qquad +\int_{\Omega}\langle f,S^{-1}_{F}F(\omega)\rangle\langle D(\omega)-S^{-1}_{F}F(\omega),f\rangle d\mu(\omega)\\
& \qquad +\int_{\Omega}\langle f,S^{-1}_{F}F(\omega)\rangle\langle S^{-1}_{F}F(\omega),f\rangle d\mu(\omega)\\
&=\int_{\Omega}\langle f,D(\omega)-S^{-1}_{F}F(\omega)\rangle\langle D(\omega)-S^{-1}_{F}F(\omega),f\rangle d\mu(\omega)\\
& \qquad +\int_{\Omega}\langle f,S^{-1}_{F}F(\omega)\rangle\langle S^{-1}_{F}F(\omega),f\rangle d\mu(\omega).
\end{align*}
Since each integral is a positive element of $C^*$-algebra $\mathcal A$, so
\begin{equation*}
\int_{\Omega}\langle f,G(\omega)\rangle\langle G(\omega),f\rangle d\mu(\omega)\leq\int_{\Omega}\langle f,D(\omega)\rangle\langle D(\omega),f\rangle d\mu(\omega),
\end{equation*}
holds for any dual $D$ of $F$ and every $f\in U$.

Conversly, for $D=S^{-1}_{F}F$, we have
\begin{equation*}
\int_{\Omega}\langle f,G(\omega)\rangle\langle G(\omega),f\rangle d\mu(\omega)\leq\int_{\Omega}\langle f,S^{-1}_{F}F(\omega)\rangle\langle S^{-1}_{F}F(\omega),f\rangle d\mu(\omega).
\end{equation*}
Also by first part of the proof,
\begin{equation*}
\int_{\Omega}\langle f,S^{-1}_{F}F(\omega)\rangle\langle S^{-1}_{F}F(\omega),f\rangle d\mu(\omega)\leq\int_{\Omega}\langle f,G(\omega)\rangle\langle G(\omega),f\rangle d\mu(\omega).
\end{equation*}
Therefore, $G=S^{-1}_{F}F$.
\end{proof}
\section{Sum of dual continuous frames}
In this section, we give the conditions under which the sum of two duals of a continuous frame is also a dual of original frame.

First, we investigate the conditions so that the sum of a continuous frame with its dual under the influence of adjointable mappings becomes a continuous frame.
\begin{theorem}\label{sum1}
Let $F:\Omega \to U$ be a continuous frame for Hilbert $C^{\ast}$-module $U$ over a unital $C^*$-algebra $\mathcal A$ with bounds $A_{F},B_{F}$ and $G:\Omega \to U$ is a dual of $F$ with bounds $A_{G},B_{G}$. Assume that $L_{1},L_{2}\in End^{*}(U)$ such that $L_{1}L_{2}^{*}= I_{U}$. Then $L_{1}F+L_{2}G$ is a continuous frame for $U$.
\end{theorem}

\begin{proof}
For every $f\in U$, we have
\begin{align*}
\int_{\Omega}\langle f,L_{1}F(\omega)+ & L_{2}G(\omega)\rangle\langle L_{1}F(\omega)+L_{2}G(\omega),f\rangle d\mu(\omega)\\ &=\int_{\Omega}\langle f,L_{1}F(\omega)\rangle\langle L_{1}F(\omega),f\rangle d\mu(\omega)+\int_{\Omega}\langle f,L_{1}F(\omega)\rangle\langle L_{2}G(\omega),f\rangle d\mu(\omega)\\
&+\int_{\Omega}\langle f,L_{2}G(\omega)\rangle\langle L_{1}F(\omega),f\rangle d\mu(\omega)+\int_{\Omega}\langle f,L_{2}G(\omega)\rangle\langle L_{2}G(\omega),f\rangle d\mu(\omega)\\
&=\int_{\Omega}\langle f,L_{1}F(\omega)\rangle\langle L_{1}F(\omega),f\rangle d\mu(\omega)+\langle L_{1}^{*}f,L_{2}^{*}f\rangle\\
&+ \langle L_{2}^{*}f,L_{1}^{*}f\rangle+\int_{\Omega}\langle f,L_{2}G(\omega)\rangle\langle L_{2}G(\omega),f\rangle d\mu(\omega)\\
&=\int_{\Omega}\langle f,L_{1}F(\omega)\rangle\langle L_{1}F(\omega),f\rangle d\mu(\omega) +2\langle f,f\rangle +\int_{\Omega}\langle f,L_{2}G(\omega)\rangle\langle L_{2}G(\omega),f\rangle d\mu(\omega).
\end{align*}
Then
\begin{equation*}
2\langle f,f\rangle\leq\int_{\Omega}\langle f,L_{1}F(\omega)+L_{2}G(\omega)\rangle\langle L_{1}F(\omega)+L_{2}G(\omega),f\rangle d\mu(\omega).
\end{equation*}
Also, since $L_{1},L_{2}$ are bounded and $\mathcal A$-linear so by Theorem \ref{BL}, there exist $M_{1},M_{2}>0$ such that
\begin{align*}
\int_{\Omega}\langle f,L_{1}F(\omega)+L_{2}G(\omega)\rangle\langle L_{1}F(\omega)+L_{2}G(\omega),f\rangle d\mu(\omega)& \leq B_{F}\langle L_{1}^{*}f,L_{1}^{*}f\rangle +2\langle f,f\rangle +B_{G}\langle L_{2}^{*}f,L_{2}^{*}f\rangle\\
& \leq(B_{F}M_{1}+2+B_{G}M_{2})\langle f,f\rangle.
\end{align*}
\end{proof}
Due to Theorem \ref{sum1}, the following results hold.
\begin{corollary}
Let $F:\Omega \to U$ be a continuous frame for Hilbert $C^{\ast}$-module $U$ over a unital $C^*$-algebra $\mathcal A$ and $G:\Omega \to U$ is a dual of $F$. Assume that $L,L_{1},L_{2}\in End^{*}(U)$ such that $L_{1}L_{2}^{*}= I_{U}$. Then the following statements hold.\\
$(1)$\; $(L_{1}F,L_{2}G)$ is a dual pair.\\
$(2)$\; If $L$ is a unitary operator, then $LF+LG$  is a continuous frame for $U$.\\
$(3)$\; $F+G$ is a continuous frame for $U$.\\
$(4)$\; $L_{1}F+L_{2}G$ is a Riesz-type frame for $U$ if and only if $T^{*}_{F}L^{*}_{1}+T^{*}_{G}L^{*}_{2}$ is surjective, where $T_{F}$and $T_{G}$ are the pre-frames operators of $F$ and $G$, respectively.
\end{corollary}
\begin{proposition}
Let $F:\Omega \to U$ be a continuous frame for Hilbert $C^{\ast}$-module $U$ over a unital $C^*$-algebra $\mathcal A$ and $G:\Omega \to U$ is a dual of $F$. If there exists $L\in End^{*}(U)$ such that $LG$ is a dual of $F$ then $L=I_{U}$.
\end{proposition}

\begin{proof}
For each $f,g\in U$ we have
\begin{align*}
\langle f,g\rangle &=\int_{\Omega}\langle f,F(\omega)\rangle\langle LG(\omega),g\rangle d\mu(\omega)\\
&= \int_{\Omega}\langle f,F(\omega)\rangle\langle G(\omega),L^{*}g\rangle d\mu(\omega)\\
&= \langle f,L^{*}g\rangle =\langle Lf,g\rangle.
\end{align*}
Hence $L=I_{U}$. 
\end{proof}
The following lemma is necessary to prove the next theorems.
\begin{lemma}\label{aF-B}
Let $F:\Omega \to U$ be a continuous Bessel mapping for Hilbert $C^{\ast}$-module $U$ over a unital $C^*$-algebra $\mathcal A$ with the bound $B$ and $a\in \mathcal A$. Then $aF:\Omega \to U$ is a continuous Bessel mapping for $U$.
\end{lemma}
\begin{proof}
\begin{align*}
\Vert\int_{\Omega}\langle f,aF(\omega)\rangle\langle aF(\omega),f\rangle d\mu(\omega)\Vert &=\Vert\lbrace\langle f,F(\omega)\rangle a^{*}\rbrace_{\omega\in\Omega}\Vert^{2}\\
&= \Vert\lbrace\langle f,F(\omega)\rangle\rbrace_{\omega\in\Omega}\;a^{*}\Vert^{2}\\
& \leq\Vert\lbrace\langle f,F(\omega)\rangle\rbrace_{\omega\in\Omega}\Vert^{2}\Vert a^{*}\Vert^{2}\leq B\Vert f\Vert^{2}\Vert a\Vert^{2},
\end{align*}
for all $f\in U$.
\end{proof}
\begin{remark}
Let $F:\Omega \to U$ be a continuous frame for Hilbert $C^{\ast}$-module $U$ over a unital $C^*$-algebra $\mathcal A$ and $a\in \mathcal A$. It is easily to see that:\\
$(1)$\; If $a$ is unitary, then $aF$ is a continuous frame for $U$ and $S_{aF}=S_{F}$.\\
$(2)$\; If $a$ is unitary and $G:\Omega \to U$ is a dual of $F$, then $aG$ is a dual of $aF$.\\
$(3)$\; If  $ab=ba$, for all $b\in \mathcal A$, then $aF$ is a continuous frame for $U$ and $S_{aF}=\vert a\vert^{2}S_{F}$ and $T_{aF}=aT_{F}$.
\end{remark}
In the following theorems, we check the conditions that the sum of two duals of a continuous frame becomes a dual of it.
\begin{theorem}
Let $F:\Omega \to U$ be a continuous frame for Hilbert $C^{\ast}$-module $U$ over a unital $C^*$-algebra $\mathcal A$  and $a_{1},a_{2}\in \mathcal A$ such that $a_{1}+a_{2}=1_{\mathcal A}$ and $a_{i}b=ba_{i}$, for $i=1,2$ and all $b\in \mathcal A$. If $G,K:\Omega \to U$ are duals of $F$, then $a_{1}G+a_{2}K$ is a dual of $F$.
\end{theorem}

\begin{proof}
By Lemma \ref{aF-B}, the mappings $a_{1}G$ and $a_{2}K$ are continuous Bessel mapping. Also
\begin{align*}
\int_{\Omega}\langle f,a_{1}G(\omega)+a_{2}K(\omega)\rangle \langle F(\omega),g\rangle d\mu(\omega) &= \int_{\Omega}\langle f,a_{1}G(\omega)\rangle\langle F(\omega),g\rangle d\mu(\omega) \\
&\;\;\;\;\;\; +\int_{\Omega}\langle f,(1_{\mathcal A}-a_{1})K(\omega)\rangle\langle F(\omega),g\rangle d\mu(\omega) \\
&= a^{*}_{1}\langle f,g\rangle +(1_{\mathcal A}-a^{*}_{1})\langle f,g\rangle =\langle f,g\rangle,
\end{align*}
for all $f,g\in U$.
\end{proof}
\begin{theorem}\label{Sum2}
Let $F:\Omega \to U$ be a continuous frame for Hilbert $C^{\ast}$-module $U$ over a unital $C^*$-algebra $\mathcal A$ with the pre-frame operator $T_{F}$. Also suppose that $G,K:\Omega \to U$ are two duals of $F$ and $V_{1},V_{2}\in End^{*}(U)$. Then $V_{1}G+V_{2}K$ is a dual of $F$ if and only if $V_{1}+V_{2}=I_{U}$.
\end{theorem}
\begin{proof}
Assume that $T_{G}$ and $T_{K}$ are the pre-frame operators of $G$ and $K$ respectively, and $F^{\prime}=V_{1}G+V_{2}K$ is a dual of $F$ with the pre-frame operator $T_{F^{\prime}}$. Then $T_{F}T^{*}_{F^{\prime}}=I_{U}$ and
\begin{align*}
T_{F}(V_{1}T_{G}+V_{2}T_{K})^{*}=I_{U} &\;\;\Longrightarrow\;\;T_{F}T^{*}_{G}V^{*}_{1}+T_{F}T^{*}_{K}V^{*}_{2}=I_{U}\\
& \;\;\Longrightarrow\;\;V^{*}_{1}+V^{*}_{2}=I_{U}\\
& \;\;\Longrightarrow\;\;V_{1}+V_{2}=I_{U}.
\end{align*}

Conversely, assume that $V_{1}+V_{2}=I_{U}$. Clearly $V_{1}G+V_{2}K$ is a continuous Bessel mapping for $U$. Then
\begin{align*}
\int_{\Omega}\langle f,V_{1}G(\omega)+V_{2}K(\omega)\rangle F(\omega)d\mu(\omega) &= \int_{\Omega}\langle f,V_{1}G(\omega)\rangle F(\omega)d\mu(\omega) \\
&\qquad +\int_{\Omega}\langle f,(I_{U}-V_{1})K(\omega)\rangle F(\omega)d\mu(\omega) \\
&= \int_{\Omega}\langle f,V_{1}G(\omega)\rangle F(\omega)d\mu(\omega)\\
&\qquad +\int_{\Omega}\langle f,K(\omega)\rangle F(\omega)d\mu(\omega) -\int_{\Omega}\langle f,V_{1}K(\omega)\rangle F(\omega)d\mu(\omega)\\
&=  V^{*}_{1}f+f-V^{*}_{1}f=f,
\end{align*}
for all $f\in U$.
\end{proof}
\begin{corollary}
Let $F:\Omega \to U$ be a continuous frame for Hilbert $C^{\ast}$-module $U$ over a unital $C^*$-algebra $\mathcal A$  and $\alpha ,\beta\in\mathbb{C}$ such that $\alpha +\beta=1$. If $G,K:\Omega \to U$ are duals of $F$, then $\alpha G+\beta K$ is a dual of $F$.
\end{corollary}

\begin{proof}
Applying Theorem \ref{Sum2}, it is enough to set $V_{1}=\alpha I_{U}$ and $V_{2}=\beta I_{U}$.
\end{proof}
\subsection*{Conflict of interest statement and data availability:}
We have no conflict of interest.

\end{document}